\documentclass[10pt,a4paper]{article}
\usepackage[T1]{fontenc}
\usepackage[utf8]{inputenc}
\usepackage{authblk}
\usepackage{amsmath,amssymb,amsthm,esint,bm}
\usepackage{mathrsfs}
\usepackage{bookmark}
\usepackage{amsmath}
\allowdisplaybreaks[3]

\newtheorem{definition}{Definition}[section]
\newtheorem{theorem}[definition]{Theorem}
\newtheorem{lemma}[definition]{Lemma}

\theoremstyle{remark}
\newtheorem{remark}[definition]{Remark}
\numberwithin{equation}{section}

\title{Radial symmetry and Liouville theorem for master equations}

\author[a,d]{Lingwei Ma}
\author[b,d]{Yahong Guo}
\author[c]{Zhenqiu Zhang\thanks{Corresponding author.}}

\affil[a]{School of Mathematical Sciences, Tianjin Normal University, Tianjin, 300387, P.~R. ~China}
\affil[b]{School of Mathematical Sciences, Nankai University, Tianjin, 300071, P.~R.~China}
\affil[c]{School of Mathematical Sciences and LPMC, Nankai University, Tianjin, 300071, P.~R.~China}
\affil[d]{Department of Mathematical Sciences, Yeshiva University, New York, NY, 10033, USA}

\usepackage{hyperref}
\begin{document}
\maketitle
\footnotetext[1]{E-mail: mlw1103@outlook.com (L. Ma), 1120200036@mail.nankai.edu.cn (Y. Guo), zqzhang@nankai.edu.cn (Z. Zhang). }

\begin{abstract}
This paper has two primary objectives. The first one is to demonstrate that the solutions of master equation
\begin{equation*}
(\partial_t-\Delta)^s u(x,t) =f(u(x, t)), \,\,(x, t)\in B_1(0)\times \mathbb{R},
\end{equation*}
subject to the vanishing exterior condition, are radially symmetric and strictly decreasing with respect to the origin in $B_1(0)$ for any $t\in \mathbb{R}$. Another one is to establish the Liouville theorem for homogeneous master equation
\begin{equation*}
    (\partial_t-\Delta)^s u(x,t)=0 ,\,\,  \mbox{in}\,\,  \mathbb{R}^n\times\mathbb{R},
\end{equation*}
which states that all bounded solutions must be constant. We propose a new methodology for a direct method of moving planes applicable to the fully fractional heat operator $(\partial_t-\Delta)^s$, and the proof of our main results based on this direct method involves the perturbation technique, limit argument as well as Fourier transform.
This study opens up a way to investigate the geometric behavior of master equations, and provides valuable insights for establishing qualitative properties of solutions and even for deriving important Liouville theorems for other types of fractional order parabolic equations.


\end{abstract}

\section{Introduction}
\label{1}

The objective of this paper is to explore the qualitative properties of solutions to space-time dual nonlocal equations involving the fractional powers of the heat operator $(\partial_t-\Delta)^s$. More specifically, we establish the radial symmetry and monotonicity of solutions for the following master equation in a unit ball
\begin{equation}\label{modelB}
\left\{\begin{array}{ll}
(\partial_t-\Delta)^s u(x,t) =f(u(x, t)), &(x, t)\in B_1(0)\times \mathbb{R}, \\
u(x, t)\equiv0, &(x, t) \in  B^c_1(0) \times \mathbb{R},
\end{array}
\right.
\end{equation}
for any $t\in \mathbb{R}$,
and prove the Liouville theorem for the homogeneous master equations
\begin{equation}\label{model}
    (\partial_t-\Delta)^s u(x,t)=0 ,\,\,  \mbox{in}\,\,  \mathbb{R}^n\times\mathbb{R}
\end{equation}
in the whole space.

The first proposal of such a fully fractional heat operator $(\partial_t-\Delta)^s$ is attributed to the mathematician Marcel Riesz, who introduced it in \cite{Riesz}. This
nonlocal operator can be defined in the following pointwise form
\begin{equation}\label{nonlocaloper}
(\partial_t-\Delta)^s u(x,t)
:=C_{n,s}\int_{-\infty}^{t}\int_{\mathbb{R}^n}
  \frac{u(x,t)-u(y,\tau)}{(t-\tau)^{\frac{n}{2}+1+s}}e^{-\frac{|x-y|^2}{4(t-\tau)}}\operatorname{d}\!y\operatorname{d}\!\tau,
\end{equation}
where $0<s<1$, the integral in $y$ is taken in the Cauchy principal value sense and the normalization positive constant $$C_{n,s}=\frac{1}{(4\pi)^{\frac{n}{2}}|\Gamma(-s)|},$$
with $\Gamma(\cdot)$ denoting the Gamma function. The singular integral in \eqref{nonlocaloper} is well defined in $\mathbb{R}^n\times\mathbb{R}$ provided
 $$u(x,t)\in C^{2s+\epsilon,s+\epsilon}_{x,\, t,\, {\rm loc}}(\mathbb{R}^n\times\mathbb{R}) \cap \mathcal{L}(\mathbb{R}^n\times\mathbb{R})$$
for some $\varepsilon>0$, where
the slowly increasing function space $\mathcal{L}(\mathbb{R}^n\times\mathbb{R})$ is defined by
$$ \mathcal{L}(\mathbb{R}^n\times\mathbb{R}):=\left\{u(x,t) \in L^1_{\rm loc} (\mathbb{R}^n\times\mathbb{R}) \mid \int_{-\infty}^t \int_{\mathbb{R}^n} \frac{|u(x,\tau)|e^{-\frac{|x|^2}{4(t-\tau)}}}{1+(t-\tau)^{\frac{n}{2}+1+s}}\operatorname{d}\!x\operatorname{d}\!\tau<\infty,\,\, \forall \,t\in\mathbb{R}\right\},$$
and the definition of the local parabolic H\"{o}lder space $C^{2s+\epsilon,s+\epsilon}_{x,\, t,\, {\rm loc}}(\mathbb{R}^n\times\mathbb{R})$ will be specified in Section \ref{2}\,. In particular, if $u$ is bounded, we can ensure the integrability of \eqref{nonlocaloper} by assuming only that $u$ is local parabolic H\"{o}lder continuous.
We notice that the operator $(\partial_t-\Delta)^s$ is nonlocal both in space and time, since the value of $(\partial_t-\Delta)^s u$ at a given point $(x,t)$ depends on the values of $u$ over the whole $\mathbb{R}^n$ and even on all the past time before $t$.
The intriguing aspect of this problem is that applying the space-time nonlocal operator $(\partial_t-\Delta)^s$ to a function that only depends on either space or time, it reduces to a familiar fractional order operator, as discussed in \cite{ST}. More precisely, if $u$ is only a function of $x$, then
 \begin{equation*}
   (\partial_t-\Delta)^s u(x)=(-\Delta)^s u(x),
 \end{equation*}
where $(-\Delta)^s$ is the well-known fractional Laplacian of order $2s$.
While if $u=u(t)$, then
 \begin{equation*}
   (\partial_t-\Delta)^s u(t)=\partial_t^s u(t),
 \end{equation*}
where $\partial_t^s$ is usually denoted by $D_{\rm left}^s$, representing the Marchaud left fractional derivative of order $s$, defined as
\begin{equation*}
D_{\rm left}^s u(t)=\frac{1}{|\Gamma(-s)|}\int_{-\infty}^t \frac{u(t)-u(\tau)}{(t-\tau)^{1+s}}\operatorname{d}\!\tau.
\end{equation*}
Moreover, it should be noted that as $s$ tends to $1$ from the left side,
the fractional power of heat operator $(\partial_t-\Delta)^s$
converges to the local heat operator $\partial_t-\Delta$ (cf. \cite{FNW}).

The master equation has a wide range of applications in physical and biological phenomena, such as anomalous diffusion \cite{KBS}, chaotic dynamics \cite{Z}, biological invasions \cite{BRR}, among others. In addition to these areas, it has also been employed in the financial field \cite{RSM}, where it can model the correlation between waiting times and price jumps in transactions. From a probabilistic perspective, the master equation plays a crucial role in the theory of continuous time random walk, where $u$ represents the distribution of particles that make random jumps simultaneously with random time lags (cf. \cite{MK}).
This is in contrast to the nonlocal parabolic equations
\begin{equation}\label{frac-para}
  \partial_t u+(-\Delta)^s u=f
\end{equation}
or the dual fractional parabolic equation
\begin{equation}\label{frac-para1}
  \partial_t^\alpha u+(-\Delta)^s u=f,
\end{equation}
where jumps are independent of the waiting times.
In other words, the master equation takes into account the strong correlation between the waiting times and the particle jumps, whereas the nonlocal parabolic equation \eqref{frac-para} or the dual fractional parabolic equation \eqref{frac-para1} do not.
It is evident that the master equation is of great importance in various fields, and continuous research on it can drive us towards a deeper understanding of complex phenomena.

Nowadays substantial progress has been made
in a series of remarkable papers \cite{ACM, CS2, ST} investigating the existence, uniqueness and regularity of solutions to master equations. The primary approach used in such studies is the extension method introduced by Caffarelli and Silvestre \cite{CS}, which extends the master equation to a local degenerate parabolic equation in a higher dimensional space.
As far as we know,
there is limited understanding of the geometric behavior of solutions to the master equation.
This lack of results can mainly be attributed to the challenges posed by the non-locality and strong correlation of the operator $(\partial_t-\Delta)^s$ within the framework of master equation.
The only related paper we are aware of is \cite{CM1}, in which Chen and Ma utilized a direct sliding method to establish that the entire solution of the master equation in \eqref{modelB} is monotone increasing and one-dimensional symmetric in $\mathbb{R}^n\times\mathbb{R}$, and thus proved the Gibbons' conjecture in this context.

The method of moving planes, which was introduced by Alexandroff in \cite{H}, is a commonly used technique to study the monotonicity and symmetry of solutions to local elliptic and parabolic equations. However, this approach cannot be directly applied to pseudo-differential equations involving the fractional Laplacian, due to the non-locality of this operator.
One effective method is to combine the aforementioned extension method, which enables us to apply the traditional method of moving planes designed for local equations to the extended problem, thereby establishing the properties of solutions.
Another useful approach is to convert
the given pseudo-differential equations into their equivalent integral equations.  By doing so, one can use the method of moving planes in integral forms and the regularity lifting to investigate the properties of
solutions (cf. \cite{CLO2, CLO}). These two effective methods have been successfully employed to investigate elliptic equations involving the fractional Laplacian. However, the above two methods can only be applied to equations involving the fractional Laplacian, and sometimes one may need to perform
cumbersome calculations and to impose additional restrictions on the problem, which may not be necessary when dealing with the fractional equations directly.

A decade later, Chen, Li, and Li \cite{CLL} made further progress by introducing a direct method of moving planes. This method removed the restrictions and greatly simplified the proof process. Since then, this effective direct method has been widely applied to establish the symmetry, monotonicity, non-existence, and even to obtain estimates in a boundary layer of solutions for various elliptic equations and systems involving the fractional Laplacian, the fully nonlinear nonlocal operators, the fractional $p$-Laplacians as well as the higher order fractional operators.
For more details, please refer to \cite{CH, CL, CLLg, CWu1, MZ1, MZ2, MZ3, ZL} and the references therein.
The direct method of moving planes is noteworthy for its ability not only to explore the symmetry, monotonicity, and non-existence of positive solutions for the fractional parabolic equations of type \eqref{frac-para} (cf. \cite{CWu2, CWNH, CWW, WuC}), but also for its recent generalization to investigate the monotonicity of solutions for the dual nonlocal parabolic equations of type \eqref{frac-para1} in a half-space (cf. \cite{CM2}).

In contrast, there has been no progress in investigating the feasibility of a direct method of moving plane for the master equation, and in exploring
how to utilize this approach to derive qualitative properties of solutions to the master equation. As we have observed that the kernel of the fully fractional heat operator $(\partial_t-\Delta)^s$ possesses a radial decreasing property, thus it is hopeful to establish the direct method of moving planes, which is exactly the research objective of this paper.
The innovation of this paper is that it overcomes the difficulties arising from the non-locality and strong correlation of the operator $(\partial_t-\Delta)^s$, and avoids the heavy reliance on extension method as in classical approach when studying the master equations. By establishing various maximum principles, we have developed a direct method of moving plane applicable to the master equation, which has allowed us to obtain the radial symmetry and strict monotonicity of solution $u(x,t)$ to the master equation \eqref{modelB} in a unit ball for any $t\in\mathbb{R}$. More surprisingly, we have also applied this direct method to establish the Liouville theorem of the homogeneous master equation in $\mathbb{R}^n\times\mathbb{R}$.

To illustrate the main results of this paper, we start by presenting the notation that will be used throughout the subsequent sections. Let $x_1$ be any given direction in $\mathbb{R}^n$, $$T_\lambda=\{x=(x_1,x')\in \mathbb R^n \mid x_1=\lambda\,\,\mbox{for}\,\,\lambda \in \mathbb R\}$$
be the moving planes perpendicular to $x_1$-axis,
$$ \Sigma_\lambda= \{x\in \mathbb R^n \mid x_1< \lambda \}\,\,\mbox{and} \,\,\Omega_\lambda= \{x\in B_1(0) \mid x_1< \lambda \}$$
be the region to the left of the hyperplane $T_\lambda$ in $\mathbb{R}^n$ and in $B_1(0)$ respectively. We denote the reflection of $x$ with respect to the hyperplane $T_\lambda$ as
$$x^\lambda=(2\lambda-x_1, x_2,\cdots, x_n).$$
Let
$u_{\lambda}(x,t)=u(x^{\lambda},t)$, we define
\begin{equation*}
w_{\lambda} (x,t) =u_{\lambda}(x,t) - u(x,t),
\end{equation*}
which represents the comparison between the values of $u(x,t)$ and $u(x^{\lambda},t)$. It is obvious that $w_{\lambda}(x,t)$ is an antisymmetric
function of $x$ with respect to the hyperplane $T_\lambda$.

We are now ready to present the main results of this paper. Our primary outcome is the narrow region principle for antisymmetric functions in bounded domains.
\begin{theorem}\label{NRP}
Let $\Omega$ be a bounded narrow domain containing in the narrow slab $\{x\in\Sigma_\lambda \mid \lambda-l<x_1<\lambda\}$  with some small $l$. Suppose that
 $$w(x,t)\in C^{2s+\epsilon,s+\epsilon}_{x,\, t,\, {\rm loc}}(\Omega\times\mathbb{R}) \cap \mathcal{L}(\mathbb{R}^n\times\mathbb{R})$$
is lower semi-continuous up to the boundary $\partial\Omega$ and bounded from below in $\Omega\times\mathbb{R}$, and satisfies
\begin{equation}\label{NRP2}
\left\{
\begin{array}{ll}
    (\partial_t-\Delta)^s w(x,t)=c(x,t) w(x,t) ,~   &(x,t) \in  \Omega\times\mathbb{R}  , \\
  w(x,t)\geq 0 , ~ &(x,t)  \in (\Sigma_\lambda  \backslash \Omega) \times\mathbb{R} , \\
  w(x,t)=- w(x^\lambda,t), &(x,t) \in \Sigma_\lambda\times\mathbb{R},
\end{array}
\right.
\end{equation}
where the coefficient function $c(x,t)$ has a uniformly upper bound $C_0$.

Then
\begin{equation}\label{NRP3}
w(x,t)\geq0, \,\,\mbox{in} \,\, \Sigma_\lambda\times\mathbb{R}
\end{equation}
for sufficiently small $l$.
 Furthermore, if $w(x,t)$ attains zero at some point $(x^0,t_0)\in\Omega\times\mathbb{R}$, then
\begin{equation}\label{NRP3-1}
 w(x,t)\equiv0, \,\, \mbox{in}\,\, \mathbb{R}^n\times(-\infty,t_0].
\end{equation}
\end{theorem}

The above narrow region principle is a crucial element in implementing the method of moving planes, as it provides a starting point to observe some geometric behavior of solutions.
Based on this maximum principle, we construct a direct method of moving planes that is suitable for the master equations, and then determine the maximum distance that the plane can be moved while maintaining the desired property. This allows us to demonstrate the following radial symmetry and strict monotonicity result.
\begin{theorem}\label{Ballsym}
Let $$u(x,t)\in C^{2s+\epsilon,s+\epsilon}_{x,\, t}(B_1(0)\times\mathbb{R}) $$ be a positive bounded solution of
\begin{equation*}
\left\{\begin{array}{ll}
(\partial_t-\Delta)^s u(x,t) =f(u(x, t)), &(x, t)\in B_1(0)\times \mathbb{R}, \\
u(x, t)\equiv0, &(x, t) \in  B^c_1(0) \times \mathbb{R}.
\end{array}
\right.
\end{equation*}
Suppose that $f\in C^1([0,+\infty))$ satisfies $f(0)\geq 0$ and $f'(0)\leq 0$.
Then $u(x,t)$ is radially symmetric and strictly decreasing about the origin in $B_1(0)$ for any $t\in \mathbb{R}$.
\end{theorem}

To establish the Liouville theorem for the homogeneous master equation in $\mathbb{R}^n\times \mathbb{R}$, which has important applications in several areas of mathematics and physics, we first establish the following maximum principle for antisymmetric functions in unbounded domains.
\begin{theorem}\label{MPUB} Let
 $$w(x,t)\in C^{2s+\epsilon,s+\epsilon}_{x,\, t,\, {\rm loc}}(\Sigma_\lambda\times\mathbb{R}) \cap \mathcal{L}(\mathbb{R}^n\times\mathbb{R})$$
be upper semi-continuous up to the boundary $T_\lambda$ and bounded from above in $\Sigma_\lambda\times\mathbb{R}$. Suppose that
\begin{equation}\label{MPUB2}
\left\{
\begin{array}{ll}
    (\partial_t-\Delta)^sw(x,t)\leq 0 ,~   & \mbox{at the points in}\,\, \Sigma_\lambda\times\mathbb{R} \,\,\mbox{where}\,\, w(x,t)>0 , \\
    w(x,t)=- w(x^\lambda,t), &\mbox{in}\,\, \Sigma_\lambda\times\mathbb{R},
\end{array}
\right.
\end{equation}
then there holds that
\begin{equation}\label{MPUB3}
w(x,t)\leq0 \,\,\mbox{in} \,\,  \Sigma_\lambda\times\mathbb{R}.
\end{equation}
\end{theorem}

\begin{remark}
Compared to the narrow domain principle (Theorem \ref{NRP}) and the maximum principle (Theorem \ref{MPUB}), although the overall proof methods of both theorems are based on the contradiction and perturbation technique, the former only has $t$ belonging to an unbounded interval, while the latter has both $x$ and $t$ belonging to unbounded domains, and $x$ does not belong to a narrow region. In the latter case, $x$ and $t$ need to be perturbed simultaneously, and it is more effective to construct a sequence of auxiliary functions in an antisymmetric form.
\end{remark}

By utilizing this maximum principle to show that the solutions are symmetric with respect to all hyperplanes in $\mathbb{R}^n$ for any $t\in\mathbb{R}$, and combining with Fourier transform, we establish the Liouville theorem for the homogeneous master equation in $\mathbb{R}^n\times \mathbb{R}$.

\begin{theorem}\label{Liouville}
Assume that $$u(x,t)\in C^{2s+\epsilon,s+\epsilon}_{x,\, t,\,{\rm loc}}(\mathbb{R}^n\times\mathbb{R}) $$ is a bounded solution of
\begin{equation*}
(\partial_t-\Delta)^s u(x,t) =0, \,\,(x, t)\in \mathbb{R}^n\times \mathbb{R},
\end{equation*}
Then $u(x,t)$ must be a constant.
\end{theorem}

\begin{remark}
As is well known, the classical Liouville theorem states that any bounded harmonic function defined on the whole space must be a constant. In fact, this boundedness condition can be relaxed to one-sided  boundedness. However, under this relaxed condition, the Liouville theorem for the caloric function satisfying heat equation
\begin{equation}\label{localheat}
  (\partial_t-\Delta)u(x,t)=0,\,\, \mbox{in}\,\, \mathbb{R}^n\times\mathbb{R}
\end{equation}
is not valid. For instance, the function $u(x,t)=e^{x+t}$ is a positive nonconstant solution of \eqref{localheat} in $\mathbb{R}^n\times \mathbb{R}$. Since the master equations recover the classical parabolic equations in limit cases, then in this sense Theorem \ref{Liouville} can be regarded as a generalization of the Liouville theorem for the local heat equation to the case of master equations involving the fractional heat operator $(\partial_t-\Delta)^s$, where the bounded condition may not be optimal but is still reasonable.
\end{remark}

\begin{remark}
In contrast to the nonlocal elliptic equation $(-\Delta)^su(x)=0$ and the nonlocal parabolic equations  $\partial_tu(x,t)+(-\Delta)^su(x,t)=0$, where the Liouville theorem can be directly proven by the maximum principles in unbounded domains. But for the space-time nonlocal equation \eqref{model}, where the operator $(\partial_t-\Delta)^s$ with respect to the time variable $t$ is also nonlocal,
such maximum principle can only determine that the solution is solely dependent on the variable $t$. To further prove that the solution must be a constant, Fourier transform is applied.
\end{remark}

To conclude this section, we will provide a brief outline of the structure of this paper. In Section \ref{2}\,, we present the definition of parabolic H\"{o}lder space, as well as some auxiliary results that are necessary for proving our main theorems.
Section \ref{3} is devoted to demonstrating two maximum principle: the narrow domain principle (Theorem \ref{NRP}) and the maximum principle in unbounded domains (Theorem \ref{MPUB}).
By utilizing these maximum principles, we establish a direct method of moving planes for the master equation, which enables us to complete the proof of Theorem \ref{Ballsym} and Theorem \ref{Liouville} in the last section.

\section{Preliminaries}\label{2}

In this section, we present the definition of the parabolic H\"{o}lder space and collect some useful preliminary estimates, which are necessary for establishing our main results. Throughout this paper, we use $C$ to denote a general constant whose value may vary from line to line.

Now we start by stating the definition of parabolic H\"{o}lder space $C^{2\alpha,\alpha}_{x,\, t}(\mathbb{R}^n\times\mathbb{R})$ (cf. \cite{Kry}) as follows.
\begin{itemize}
\item[(i)]
When $0<\alpha\leq\frac{1}{2}$, if
$u(x,t)\in C^{2\alpha,\alpha}_{x,\, t}(\mathbb{R}^n\times\mathbb{R})$, then there exists a constant $C>0$ such that
\begin{equation*}
  |u(x,t)-u(y,\tau)|\leq C\left(|x-y|+|t-\tau|^{\frac{1}{2}}\right)^{2\alpha}
\end{equation*}
for any $x,\,y\in\mathbb{R}^n$ and $t,\,\tau\in \mathbb{R}$.
\item[(ii)]
When $\frac{1}{2}<\alpha\leq1$, we say that
$$u(x,t)\in C^{2\alpha,\alpha}_{x,\, t}(\mathbb{R}^n\times\mathbb{R}):=C^{1+(2\alpha-1),\alpha}_{x,\, t}(\mathbb{R}^n\times\mathbb{R}),$$ if $u$ is $\alpha$-H\"{o}lder continuous in $t$ uniformly with respect to $x$ and its gradient $\nabla_xu$ is $(2\alpha-1)$-H\"{o}lder continuous in $x$ uniformly with respect to $t$ and $(\alpha-\frac{1}{2})$-H\"{o}lder continuous in $t$ uniformly with respect to $x$.
\item[(iii)] While for $\alpha>1$, if
$u(x,t)\in C^{2\alpha,\alpha}_{x,\, t}(\mathbb{R}^n\times\mathbb{R}),$
then it means that
$$\partial_tu,\, D^2_xu \in C^{2\alpha-2,\alpha-1}_{x,\, t}(\mathbb{R}^n\times\mathbb{R}).$$
\end{itemize}
In addition, we can analogously define the local parabolic H\"{o}lder space $C^{2\alpha,\alpha}_{x,\, t,\, \rm{loc}}(\mathbb{R}^n\times\mathbb{R})$.

In the following, we provide the boundedness estimates for the nonlocal operators $\partial_t^s$ and $(\partial_t-\Delta)^s$ acting on smooth functions, respectively. These estimates are repeatedly employed in establishing our main results.
\begin{lemma}{\rm{(cf. \cite{CM2})}}\label{lem1}
Let $$\eta_r(t)\in C_0^\infty\left((-r^{2}+t_0,r^{2}+t_0)\right),$$
For some $t_0\in \mathbb{R}$ and $r>0$,  then there exists a positive constant $C$ that depends only on $s$ such that
$$|\partial_t^s\eta_r(t)|\leq \frac{C}{r^{2s}} \,\, \mbox{in}\,\, (-r^{2}+t_0,r^{2}+t_0).$$
\end{lemma}

\begin{lemma}{\rm{(cf. \cite{CM1})}}\label{lem2}
Let $$\phi_r(x,t)\in C_0^\infty\left(B_r(x^0)\times(-r^{2}+t_0,r^{2}+t_0)\right)$$
for some $(x^0,t_0)\in \mathbb{R}^n \times \mathbb{R}$ and $r>0$, then there exists a positive constant $C$ that depends only on $s$ and $n$ such that
$$|(\partial_t-\Delta)^s\phi_r(x,t)|\leq \frac{C}{r^{2s}} \,\, \mbox{in}\,\, B_r(x^0)\times(-r^{2}+t_0,r^{2}+t_0).$$
\end{lemma}

\section{Various maximum principles for antisymmetric functions}\label{3}

In this section, we demonstrate various maximum principles for antisymmetric functions, namely, Theorem \ref{NRP} and Theorem \ref{MPUB}\,. It will be shown in the later section that these principles play a crucial role in establishing a direct method of moving planes for the master equations.

\subsection{Narrow region principle in bounded domains}\label{3.1}

Our first objective is to prove Theorem \ref{NRP}\,, which establishes a narrow region principle for antisymmetric functions in bounded domains. This principle is a key ingredient in demonstrating the radial symmetry and monotonicity of solutions for the master equation \eqref{modelB}.

\begin{proof}
[\bf Proof of Theorem \ref{NRP}\,.] \,

We first argue by contradiction to derive \eqref{NRP3}. Due to the narrow region $\Omega$
is bounded and the function $w$ is lower semi-continuous up to the boundary $\partial\Omega$, there must exist $x(t)\in \Omega$ such that
\begin{equation}\label{NRP4}
  w(x(t),t)=\displaystyle\min_{x\in \Omega}w(x,t)
\end{equation}
for each fixed $t\in\mathbb{R}$.
If \eqref{NRP3} is not valid, we may assume on the contrary that
there exists a positive constant $m$ such that
\begin{equation}\label{NRP5}
  \inf_{\Omega\times\mathbb{R}}w(x,t)=\inf_{\mathbb{R}}w(x(t),t)=-m<0.
\end{equation}
Since $t$ belongs to $\mathbb{R}$ that is an unbounded domain, the infimum of $w(x(t),t)$ with respect to $t$ may not be attainable, but there certainly exists
a sequence $\{(x(t_k),t_k)\}\subset\Omega\times\mathbb{R}$ such that
\begin{equation*}
 w(x(t_k),t_k)=-m_k\rightarrow -m \,\,\mbox{as}\,\, k\rightarrow\infty.
\end{equation*}
Let $\varepsilon_k=m-m_k$, then it is evident that $\varepsilon_k$ is nonnegative and $\varepsilon_k\rightarrow 0$ as $k\rightarrow0$. In order to remedy scenario that the infimum of $w$ with respect to $t$ may not be attained, we need to perturb $w$ with respect to the variable $t$ such that the perturbed function can attain its infimum. For this purpose,
we introduce the following auxiliary function
\begin{equation*}
  v_k(x,t)=w(x,t)-\varepsilon_k\eta_k(t),
\end{equation*}
where
\begin{equation*}
  \eta_k(t)=\eta(t-t_k)\in C_0^\infty\left((-1+t_k,1+t_k)\right)
\end{equation*}
is a smooth cut-off function that satisfies
\begin{equation}\label{NRP5-1}
\eta(t)=\left\{
\begin{array}{ll}
1  & t\in (-\frac{1}{2},\frac{1}{2}), \\
  0 , ~ &t\not\in (-1,1),
\end{array}
\right.
\end{equation}
and $0\leq\eta(t)\leq 1$.
One one hand, we have
\begin{equation*}
  v_k(x(t_k),t_k)=w(x(t_k),t_k)-\varepsilon_k=-m_k-m+m_k=-m.
\end{equation*}
On the other hand, if $(x,t)\in\Omega\times\left(\mathbb{R}\setminus (-1+t_k,1+t_k)\right)$, then it follows from \eqref{NRP5} that
\begin{equation*}
  v_k(x,t)=w(x,t)\geq-m.
\end{equation*}
Based on the above analysis and the exterior condition in \eqref{NRP2} satisfied by $w$, there exists $$(\bar{x}^k,\bar{t}_k)\in\Omega\times(-1+t_k,1+t_k)$$
 such that
\begin{equation}\label{NRP6}
 -m-\varepsilon_k\leq v_k(\bar{x}^k,\bar{t}_k)= \inf_{\Omega\times\mathbb{R}}v_k(x,t)= \inf_{\Sigma_\lambda\times\mathbb{R}}v_k(x,t)\leq-m.
\end{equation}
From this, it is not difficult to further verify that
\begin{equation}\label{NRP6-1}
  -m\leq w(\bar{x}^k,\bar{t}_k)\leq-m+\varepsilon_k=-m_k<0.
\end{equation}

Next, we derive a contradiction at the minimum point $(\bar{x}^k,\bar{t}_k)$ of $v_k$ in $\Sigma_\lambda\times\mathbb{R}$. On the one hand, by performing a direct calculation and combining the anti-symmetry of $w$ in $x$ with the decreasing property of the kernel of operator $(\partial_t-\Delta)^s$ due to $|\bar{x}^k-y|<|\bar{x}^k-y^\lambda|$, we have
\begin{eqnarray}\label{NRP7}
&&(\partial_t-\Delta)^s v_k(\bar{x}^k,\bar{t}_k)\nonumber\\
 &=&C_{n,s}\int_{-\infty}^{\bar{t}_k}\int_{\mathbb{R}^n}
  \frac{v_k(\bar{x}^k,\bar{t}_k)-v_k(y,\tau)}{(\bar{t}_k-\tau)^{\frac{n}{2}+1+s}}e^{-\frac{|\bar{x}^k-y|^2}{4(\bar{t}_k-\tau)}}\operatorname{d}\!y\operatorname{d}\!\tau \nonumber\\
  &=&C_{n,s}\int_{-\infty}^{\bar{t}_k}\int_{\Sigma_\lambda}
  \frac{v_k(\bar{x}^k,\bar{t}_k)-v_k(y,\tau)}{(\bar{t}_k-\tau)^{\frac{n}{2}+1+s}}e^{-\frac{|\bar{x}^k-y|^2}{4(\bar{t}_k-\tau)}}\operatorname{d}\!y\operatorname{d}\!\tau \nonumber\\
  && +C_{n,s}\int_{-\infty}^{\bar{t}_k}\int_{\Sigma_\lambda}
  \frac{v_k(\bar{x}^k,\bar{t}_k)-v_k(y^\lambda,\tau)}{(\bar{t}_k-\tau)^{\frac{n}{2}+1+s}}e^{-\frac{|\bar{x}^k-y^\lambda|^2}{4(\bar{t}_k-\tau)}}\operatorname{d}\!y\operatorname{d}\!\tau
 \nonumber\\
  &\leq&C_{n,s}\int_{-\infty}^{\bar{t}_k}\int_{\Sigma_\lambda}
  \frac{v_k(\bar{x}^k,\bar{t}_k)-v_k(y,\tau)}{(\bar{t}_k-\tau)^{\frac{n}{2}+1+s}}e^{-\frac{|\bar{x}^k-y^\lambda|^2}{4(\bar{t}_k-\tau)}}\operatorname{d}\!y\operatorname{d}\!\tau \nonumber\\
  && +C_{n,s}\int_{-\infty}^{\bar{t}_k}\int_{\Sigma_\lambda}
  \frac{v_k(\bar{x}^k,\bar{t}_k)-v_k(y^\lambda,\tau)}{(\bar{t}_k-\tau)^{\frac{n}{2}+1+s}}e^{-\frac{|\bar{x}^k-y^\lambda|^2}{4(\bar{t}_k-\tau)}}\operatorname{d}\!y\operatorname{d}\!\tau
 \nonumber\\
 &=&C_{n,s}\int_{-\infty}^{\bar{t}_k}\int_{\Sigma_\lambda}
  \frac{v_k(\bar{x}^k,\bar{t}_k)-w(y,\tau)+\varepsilon_k\eta_k(\tau)}{(\bar{t}_k-\tau)^{\frac{n}{2}+1+s}}e^{-\frac{|\bar{x}^k-y^\lambda|^2}{4(\bar{t}_k-\tau)}}\operatorname{d}\!y\operatorname{d}\!\tau \nonumber\\
  && +C_{n,s}\int_{-\infty}^{\bar{t}_k}\int_{\Sigma_\lambda}
  \frac{v_k(\bar{x}^k,\bar{t}_k)-w(y^\lambda,\tau)+\varepsilon_k\eta_k(\tau)}{(\bar{t}_k-\tau)^{\frac{n}{2}+1+s}}e^{-\frac{|\bar{x}^k-y^\lambda|^2}{4(\bar{t}_k-\tau)}}\operatorname{d}\!y\operatorname{d}\!\tau
 \nonumber\\
  &=&2C_{n,s}v_k(\bar{x}^k,\bar{t}_k)\int_{-\infty}^{\bar{t}_k}\int_{\Sigma_\lambda}
  \frac{e^{-\frac{|\bar{x}^k-y^\lambda|^2}{4(\bar{t}_k-\tau)}}}{(\bar{t}_k-\tau)^{\frac{n}{2}+1+s}}\operatorname{d}\!y\operatorname{d}\!\tau +2 C_{n,s}\varepsilon_k\int_{-\infty}^{\bar{t}_k}\int_{\Sigma_\lambda}
  \frac{\eta_k(\tau)e^{-\frac{|\bar{x}^k-y^\lambda|^2}{4(\bar{t}_k-\tau)}}}{(\bar{t}_k-\tau)^{\frac{n}{2}+1+s}}\operatorname{d}\!y\operatorname{d}\!\tau. \nonumber\\
\end{eqnarray}
In order to estimate two integrals  in the last line of \eqref{NRP7}, we make a change of variables
$$t=\frac{|\bar{x}^k-y^\lambda|^2}{4(\bar{t}_k-\tau)}$$
 and derive
\begin{eqnarray*}
  && C_{n,s}\int_{-\infty}^{\bar{t}_k}\int_{\Sigma_\lambda}
  \frac{e^{-\frac{|\bar{x}^k-y^\lambda|^2}{4(\bar{t}_k-\tau)}}}{(\bar{t}_k-\tau)^{\frac{n}{2}+1+s}}\operatorname{d}\!y\operatorname{d}\!\tau  \\
   &=& \frac{1}{(4\pi)^{\frac{n}{2}}|\Gamma(-s)|} \int_{\Sigma_\lambda}\int_{0}^{+\infty}\frac{e^{-t}}{(\frac{|\bar{x}^k-y^\lambda|^2}{4t})^{\frac{n}{2}+1+s}}\frac{|\bar{x}_k-y^\lambda|^2}{4t^2}\operatorname{d}\!t \operatorname{d}\!y\\
      &=& \frac{4^s}{\pi^{\frac{n}{2}}|\Gamma(-s)|} \int_{\Sigma_\lambda}\frac{1}{|\bar{x}^k-y^\lambda|^{n+2s}}\int_{0}^{+\infty}e^{-t}t^{\frac{n}{2}+s-1}\operatorname{d}\!t \operatorname{d}\!y \\
      &=&\frac{4^s\Gamma(\frac{n}{2}+s)}{\pi^{\frac{n}{2}}|\Gamma(-s)|} \int_{\Sigma_\lambda}\frac{1}{|\bar{x}^k-y^\lambda|^{n+2s}} \operatorname{d}\!y.
\end{eqnarray*}
Substituting the above equality into \eqref{NRP7}, and applying $\eta_k\in[0,1]$ and \eqref{NRP6} to arrive at
\begin{eqnarray}\label{NRP7-1}
  &&(\partial_t-\Delta)^s v_k(\bar{x}^k,\bar{t}_k)\nonumber\\
   &\leq&Cv_k(\bar{x}^k,\bar{t}_k)\int_{\Sigma_\lambda}\frac{1}{|\bar{x}^k-y^\lambda|^{n+2s}} \operatorname{d}\!y + C\varepsilon_k\int_{\Sigma_\lambda}\frac{1}{|\bar{x}^k-y^\lambda|^{n+2s}} \operatorname{d}\!y \nonumber\\
    &\leq&\frac{Cv_k(\bar{x}^k,\bar{t}_k)}{l^{2s}}+\frac{C\varepsilon_k}{l^{2s}}\nonumber\\
    &\leq&-\frac{Cm}{l^{2s}}+\frac{C\varepsilon_k}{l^{2s}}.
\end{eqnarray}

On the other hand, starting from the differential equation in \eqref{NRP2} and combining the fact that $C(x,t)$ has a uniformly upper bound $C_0$, \eqref{NRP6-1} with Lemma \ref{lem1}\,, we obtain
\begin{eqnarray}\label{NRP7-2}
  (\partial_t-\Delta)^s v_k(\bar{x}^k,\bar{t}_k)&=& (\partial_t-\Delta)^s w(\bar{x}^k,\bar{t}_k)-\varepsilon_k(\partial_t-\Delta)^s \eta(\bar{t}_k) \nonumber\\
  &=& c(\bar{x}^k,\bar{t}_k)w(\bar{x}^k,\bar{t}_k)-\varepsilon_k\partial_t^s\eta(\bar{t}_k)\nonumber\\
  &\geq&-C_0m-C\varepsilon_k.
\end{eqnarray}
Then a combination of \eqref{NRP7-1} and \eqref{NRP7-2} yields that
\begin{equation*}
  -C_0m\leq-\frac{Cm}{l^{2s}}+\frac{C\varepsilon_k}{l^{2s}}+C\varepsilon_k\rightarrow-\frac{Cm}{l^{2s}},
\end{equation*}
as $k\rightarrow\infty$.
Dividing both side of the preceding inequality by $-m$, we deduce that
\begin{equation*}
  C_0\geq\frac{C}{l^{2s}},
\end{equation*}
which is a contradiction for sufficiently small $l$. Hence, we conclude that \eqref{NRP3} is true.

In the sequel, we remain to demonstrate the validity of \eqref{NRP3-1}. It follows from \eqref{NRP3} that
$$w(x^0,t_0)=\min_{\Sigma_\lambda\times\mathbb{R}}w(x,t)=0.$$
The equation in \eqref{NRP2} obviously implies that
\begin{equation}\label{NRP9}
  (\partial_t-\Delta)^s w(x^0,t_0)=0.
\end{equation}
Besides, through a straightforward calculation, we derive
\begin{eqnarray*}
  (\partial_t-\Delta)^s w(x^0,t_0) &=&-C_{n,s}\int_{-\infty}^{t_0}\int_{\mathbb{R}^n}
  \frac{w(y,\tau)}{(t_0-\tau)^{\frac{n}{2}+1+s}}e^{-\frac{|x^0-y|^2}{4(t_0-\tau)}}\operatorname{d}\!y\operatorname{d}\!\tau \\
   &=&  C_{n,s}\int_{-\infty}^{t_0}\int_{\Sigma_{\lambda}}
  \frac{w(y,\tau)}{(t_0-\tau)^{\frac{n}{2}+1+s}}\left[e^{-\frac{|x^0-y^\lambda|^2}{4(t_0-\tau)}}-e^{-\frac{|x^0-y|^2}{4(t_0-\tau)}}\right]\operatorname{d}\!y\operatorname{d}\!\tau.
\end{eqnarray*}
Since $w(x,t)\geq 0$ in $\Sigma_\lambda\times\mathbb{R}$ and
$$e^{-\frac{|x^0-y^\lambda|^2}{4(t_0-\tau)}}-e^{-\frac{|x^0-y|^2}{4(t_0-\tau)}}<0,$$
then it follows from \eqref{NRP9} that
\begin{equation*}
  w(x,t)\equiv 0 \,\, \mbox{for}\,\, (x,t)\in\Sigma_\lambda\times(-\infty,t_0].
\end{equation*}
Finally, the antisymmetry of $w(x,t)$ with respect to $x$ infers that
$$w(x,t)\equiv0 \,\, \mbox{in}\,\, \mathbb{R}^n\times(-\infty,t_0].$$
Therefore, the proof of Theorem \ref{NRP} is completed.
\end{proof}

\subsection{Maximum principle in unbounded domains}\label{3.2}

We proceed to prove Theorem \ref{MPUB}\,, which establishes the maximum principle for antisymmetric functions in unbounded domains. This principle is a fundamental ingredient in proving the Liouville theorem for homogeneous master equations.

\begin{proof}[\bf Proof of Theorem \ref{MPUB}\,.] \,
The proof goes by contradiction, if \eqref{MPUB3} is violated, since $w(x,t)$ has an upper bound in $\Sigma_\lambda\times\mathbb{R}$, then there exists a positive constant $A$ such that
\begin{equation}\label{MPUB4}
  \sup_{\Sigma_\lambda\times\mathbb{R}}w(x,t)=A>0.
\end{equation}
Note that the set $\Sigma_\lambda\times\mathbb{R}$ is unbounded, then the supremum of $w(x,t)$ may not be attained.
Even so, \eqref{MPUB4} implies that there exists a sequence $\{(x^k,t_k)\}\subset \Sigma_\lambda\times\mathbb{R}$ such that
\begin{equation}\label{MPUB4-1}
  0<w(x^k,t_k)=A_k\rightarrow A,\,\, \mbox{as}\,\, k\rightarrow\infty.
\end{equation}
Let $\varepsilon_k=A-A_k$, then the sequence $\{\varepsilon_k\}$ is nonnegative and
tends to zero as $k\rightarrow\infty$.

Given that the supremum of $w(x,t)$ with respect to both $x$ and $t$ may not be attained,  it is necessary to perturb the function $w$ with respect to both $x$ and $t$ such that the perturbed function not only attains the supremum, but also preserves the antisymmetry in $x$.
With this aim, we need to introduce an antisymmetric auxiliary function. Let
$\phi(x)\in C_{0}^{\infty}\left(\mathbb{R}^n\right)$ satisfy
\begin{equation*}
\phi(x)=
\left\{\begin{array}{r@{\ \ }c@{\ \ }ll}
&&e^{1+\frac{1}{|x|^{2}-1}}, & \ \ x\in B_{1}(0) \,, \\[0.05cm]
&&0, & \ \ x\not\in B_{1}(0)\,.
\end{array}\right.
\end{equation*}
We denote
\begin{equation*}
 \phi_{\lambda}(x)=\phi\left(x^{\lambda}\right)=
\left\{\begin{array}{r@{\ \ }c@{\ \ }ll}
&& e^{1+\frac{1}{|x^\lambda|^{2}-1}}, & \ \ x\in B_{1}\big(0^\lambda\big) \,, \\[0.05cm]
&&0, & \ \ x\not\in B_{1}\big(0^\lambda\big)\,,
\end{array}\right.
\end{equation*}
and $r_k=\operatorname{dist}(x^k,T_\lambda)$,
then it is not difficult to verify that
$$\Phi_k(x)=\phi\left(\frac{2(x-x^k)}{r_k}\right)-\phi_{\lambda}\left(\frac{2(x^\lambda-x^k)}{r_k}\right)\in C_{0}^{\infty}\left(B_{\frac{r_k}{2}}(x^k)\cup B_{\frac{r_k}{2}}\big((x^k)^\lambda\big) \right)$$
is an antisymmetric function with respect to the plane $T_{\lambda}$, and
$$\max_{x\in \mathbb{R}^n}\Phi_k(x)=\Phi_k(x^k)=\phi(x^k)=1.$$
We further select a smooth cut-off function of $t$ that
$$\eta_k(t)=\eta\left(\frac{t-t_k}{(\frac{r_k}{2})^2}\right)\in C_0^\infty\left((t_k-(\frac{r_k}{2})^2,t_k+(\frac{r_k}{2})^2)\right),$$
where $\eta$ is defined in \eqref{NRP5-1}.

Now we choose the antisymmetric auxiliary function as follows
$$V_k(x,t)=w(x,t)+\varepsilon_k\Phi_k(x)\eta_k(t).$$
Let
$$Q_{\frac{r_k}{2}}(x^k,t_k)=B_{\frac{r_k}{2}}(x^k)\times (t_k-(\frac{r_k}{2})^2,t_k+(\frac{r_k}{2})^2)$$
be a parabolic cylinder, then a straightforward calculation implies that
$$V_k(x^k,t_k)=w(x^k,t_k)+\varepsilon_k=A_k+A-A_k=A,$$
and
$$V_k(x,t)=w(x,t)\leq A\,\, \mbox{for}\,\, (x,t)\in (\Sigma_\lambda\times \mathbb{R})\setminus Q_{\frac{r_k}{2}}(x^k,t_k).$$
Thus, the auxiliary function $V_k(x,t)$ must attain its maximum value in $\Sigma_\lambda\times \mathbb{R}$. More precisely, there exists a point $(\bar{x}^k,\bar{t}_k)\in Q_{\frac{r_k}{2}}(x^k,t_k)$ such that
\begin{equation}\label{MPUB5}
 A+\varepsilon_k \geq V_k(\bar{x}^k,\bar{t}_k)=\sup_{\Sigma_\lambda\times\mathbb{R}}V_k(x,t)\geq A>0.
\end{equation}
Meanwhile, it follows from the definition of $V_k$ that
$$A\geq w(\bar{x}^k,\bar{t}_k)\geq A-\varepsilon_k=A_k>0,$$
then we can apply the differential inequality in \eqref{MPUB2} to $w$ at point $(\bar{x}^k,\bar{t}_k)$.

Next, we intend to derive a contradiction at the maximum point $(\bar{x}^k,\bar{t}_k)$ of $V_k$ in $\Sigma_\lambda\times\mathbb{R}$. On one hand,
combining the definition of $V_k$ with the differential inequality in \eqref{MPUB2} and Lemma \ref{lem2}\,, we obtain
\begin{equation}\label{MPUB6}
 (\partial_t-\Delta)^sV_k(\bar{x}^k,\bar{t}_k)=(\partial_t-\Delta)^sw(\bar{x}^k,\bar{t}_k)
 +\varepsilon_k(\partial_t-\Delta)^s\left(\Phi_k(\bar{x}^k)\eta_k(\bar{t}_k)\right)\leq\frac{C\varepsilon_k}{r_k^{2s}}.
\end{equation}
On the other hand, starting from the definition of operator $(\partial_t-\Delta)^s$ and utilizing the antisymmetry of auxiliary function $V_k$, the radial decrease of the kernel as well as \eqref{MPUB5}, we compute
\begin{eqnarray}\label{MPUB7}
  &&(\partial_t-\Delta)^s V_k(\bar{x}^k,\bar{t}_k)\nonumber\\
  &=&C_{n,s}\int_{-\infty}^{\bar{t}_k}\int_{\Sigma_\lambda}
  \frac{V_k(\bar{x}^k,\bar{t}_k)-V_k(y,\tau)}{(\bar{t}_k-\tau)^{\frac{n}{2}+1+s}}e^{-\frac{|\bar{x}_k-y|^2}{4(\bar{t}_k-\tau)}}\operatorname{d}\!y\operatorname{d}\!\tau \nonumber\\
  && +C_{n,s}\int_{-\infty}^{\bar{t}_k}\int_{\Sigma_\lambda}
  \frac{V_k(\bar{x}^k,\bar{t}_k)+V_k(y,\tau)}{(\bar{t}_k-\tau)^{\frac{n}{2}+1+s}}e^{-\frac{|\bar{x}_k-y^\lambda|^2}{4(\bar{t}_k-\tau)}}\operatorname{d}\!y\operatorname{d}\!\tau
 \nonumber\\
  &\geq&2C_{n,s}\int_{-\infty}^{\bar{t}_k}\int_{\Sigma_\lambda}
  \frac{V_k(\bar{x}^k,\bar{t}_k)}{(\bar{t}_k-\tau)^{\frac{n}{2}+1+s}}e^{-\frac{|\bar{x}_k-y^\lambda|^2}{4(\bar{t}_k-\tau)}}\operatorname{d}\!y\operatorname{d}\!\tau \nonumber\\
  &=&2\frac{4^s\Gamma(\frac{n}{2}+s)}{\pi^{\frac{n}{2}}|\Gamma(-s)|}V_k(\bar{x}^k,\bar{t}_k)\int_{\Sigma_\lambda}
  \frac{1}{|\bar{x}_k-y^\lambda|^{n+2s}}\operatorname{d}\!y \nonumber\\
  &\geq&\frac{CA}{r_k^{2s}}.
\end{eqnarray}

Finally, a combination of \eqref{MPUB6} and \eqref{MPUB7} yields that
\begin{equation*}
  A\leq C\epsilon_k,
\end{equation*}
which leads to a contradiction for sufficiently large $k$.
Hence, we conclude that \eqref{MPUB3} is valid, and then the proof of Theorem \ref{MPUB} is completed.
\end{proof}

\section{The proof of main results}\label{4}

In this section, we will show how the maximum principles established previously can be used to develop a direct method of moving planes applicable to the master equations. By means of this direct method, we complete the proof of our main results (i.e., Theorem \ref{Ballsym} and Theorem \ref{Liouville}).

\subsection{Radial symmetry of solutions in a unit ball}\label{4.1}

In this subsection, we apply the narrow region principle (Theorem \ref{NRP}) to initiate the moving plane and combine the perturbation technique with the limit argument to prove that solutions of the master equation subject to the vanishing exterior condition are radially symmetric and strictly decreasing with respect to the origin in $B_1(0)$ for any $t\in \mathbb{R}$, under appropriate assumptions on the nonhomogeneous term $f$.

\begin{proof}[\bf Proof of Theorem \ref{Ballsym}\,.] \,
To carrying out the method of moving planes, we choose $x_1$ to be any direction and let $T_\lambda$, $\Sigma_\lambda$, $\Omega_\lambda$, $x^\lambda$, and $w_\lambda$ be defined as in Section \ref{1}\,.
By a direct calculation, we have
\begin{equation}\label{MR2}
\left\{
\begin{array}{ll}
    (\partial_t-\Delta)^s w_\lambda(x,t)=C_\lambda(x,t)w_\lambda(x,t),~   &(x,t) \in  \Omega_\lambda\times\mathbb{R}  , \\
  w_\lambda(x,t)\geq 0 , ~ &(x,t)  \in (\Sigma_\lambda\backslash\Omega_\lambda)\times\mathbb{R},\\
   w_\lambda(x,t)=- w_\lambda(x^\lambda,t), ~ &(x,t)  \in \Sigma_\lambda\times\mathbb{R},\\
\end{array}
\right.
\end{equation}
where the coefficient function
\begin{equation*}
  C_\lambda(x,t)=\frac{f(u_\lambda(x,t))-f(u(x,t))}{u_\lambda(x,t)-u(x,t)}
\end{equation*}
is bounded ensured by $f\in C^1([0,+\infty))$ and the boundedness of $u$. Now we divide the proof into two steps.

\noindent \textup{\textbf{Step 1.}} Start moving the plane $T_{\lambda}$ from $x_1=-1$ to the right along $x_{1}$-axis. If $\lambda$ is sufficiently close to $-1$, then $\Omega_\lambda$ is a narrow region. Furthermore, the assumptions in Theorem \ref{Ballsym} guarantee that we can apply the narrow region principle, as established in Theorem \ref{NRP}\,, to problem \eqref{MR2}. This allows us to deduce that
\begin{equation}\label{MR3}
  w_{\lambda}(x,t)\geq0 \,\, \mbox{in}\,\, \Omega_{\lambda}\times\mathbb{R}.
\end{equation}
Note that inequality \eqref{MR3} provides a starting point to move the plane $T_\lambda$.

\noindent \textup{\textbf{Step 2.}} In the second step, we continue to move the plane $T_\lambda$ to the right along $x_{1}$-axis as long as \eqref{MR3} is valid to its limiting position. Let
 \begin{equation*}
 \lambda_0 = \sup \left\{ \lambda<0 \mid w_\mu (x,t) \geq 0,\,\, (x,t) \in \Sigma_\mu\times\mathbb{R}\,\,\mbox{for any} \,\, \mu \leq \lambda \right\},
 \end{equation*}
Our purpose is to show that
\begin{equation}\label{MR4}
 \lambda_0=0
\end{equation}
 by a contradiction argument.
Otherwise, if $\lambda_0<0$, then the definition of $\lambda_0$ and $w_{\lambda_k}(x,t)\geq 0$ in $(\Sigma_{\lambda_k}\setminus \Omega_{\lambda_k})\times\mathbb{R}$ implies that there exists a
sequence of negative numbers $\{\lambda_k\}$ with $\{\lambda_k\}\searrow \lambda_0$ such that
$$\inf_{\Sigma_{\lambda_k}\times \mathbb{R}}w_{\lambda_k}(x,t)=\inf_{\Omega_{\lambda_k}\times \mathbb{R}}w_{\lambda_k}(x,t)=-m_k<0,$$
and
$$m_k\rightarrow 0,\,\, \mbox{as}\,\, k\rightarrow\infty.$$
Since $ \mathbb{R}$ is an unbounded interval, the infimum of $w_{\lambda_k}$ with respect to $t$ may not be attained, but
there must exist a sequence $\{(x^k,t_k)\}\subset\Omega_{\lambda_k}\times\mathbb{R}$ and a nonnegative sequence $\{\varepsilon_k\}\searrow 0$ as $k\rightarrow\infty$, such that
\begin{eqnarray*}\label{MR5}
w_{\lambda_k}(x^k, t_k)=\inf_{x\in\Sigma_{\lambda_k}}w_{\lambda_k}(x, t_ k)\leq-m_k+\varepsilon_km_k<0.
\end{eqnarray*}
To address the situation where the infimum of $w_{\lambda_k}$ with respect to $t$ may not be attained, we need to introduce the following auxiliary function
\begin{eqnarray*}
 W_{k} (x, t)=w_{\lambda_k}(x,t)-\varepsilon_k m_k\eta_k(t),
\end{eqnarray*}
where
$\eta_k(t)=\eta(t-t_k)$
 and $\eta(t)$ is a smooth cut-off function defined in \eqref{NRP5-1}.
A straightforward computation leads to
$$W_{k}(x^k, t_k)\leq-m_k,$$
and $$W_k(x,t)=w_{\lambda_k}(x,t)\geqslant-m_k$$ for $|t-t_k|\geqslant1$.
Hence, the auxiliary function $W_{k}(x, t)$ could attain its minimum at some point $$(\bar{x}^k, \bar t_k)\in \Omega_{\lambda_k}\times(t_k-1,t_k+1),$$ such that
\begin{eqnarray}\label{MR6}
-m_k-\varepsilon_km_k \leq W_{k} (\bar {x}^k, \bar {t}_k)= \inf_{\Sigma_{\lambda_k}\times \mathbb{R}}
 W_{k}(x, t)\leq-m_k.
\end{eqnarray}
It follows that
\begin{equation}\label{MR7}
  -m_k\leq w_{\lambda_k}(\bar{x}^k,\bar{t}_k)\leq-m_k+\varepsilon_km_k<0.
\end{equation}

Next, we focus on the estimate of $W_k$ at its minimum point $(\bar {x}^k, \bar {t}_k)$. On one hand, analogous to the estimate of \eqref{NRP7-1}, applying the antisymmetry of $w_{\lambda_k}$ in $x$ and $|\bar{x}^k-y|<|\bar{x}^k-y^{\lambda_k}|$, we derive
\begin{eqnarray}\label{MR8}
&&(\partial_t-\Delta)^s W_k(\bar{x}^k,\bar{t}_k)\nonumber\\
  &\leq&2C_{n,s}W_k(\bar{x}^k,\bar{t}_k)\int_{-\infty}^{\bar{t}_k}\int_{\Sigma_{\lambda_k}}
  \frac{e^{-\frac{|\bar{x}^k-y^{\lambda_k}|^2}{4(\bar{t}_k-\tau)}}}{(\bar{t}_k-\tau)^{\frac{n}{2}+1+s}}\operatorname{d}\!y\operatorname{d}\!\tau +2 C_{n,s}\varepsilon_km_k\int_{-\infty}^{\bar{t}_k}\int_{\Sigma_{\lambda_k}}
  \frac{e^{-\frac{|\bar{x}^k-y^{\lambda_k}|^2}{4(\bar{t}_k-\tau)}}}{(\bar{t}_k-\tau)^{\frac{n}{2}+1+s}}\operatorname{d}\!y\operatorname{d}\!\tau \nonumber\\
   &=&CW_k(\bar{x}^k,\bar{t}_k)\int_{\Sigma_{\lambda_k}}\frac{1}{|\bar{x}^k-y^{\lambda_k}|^{n+2s}} \operatorname{d}\!y + C\varepsilon_km_k\int_{\Sigma_{\lambda_k}}\frac{1}{|\bar{x}^k-y^{\lambda_k}|^{n+2s}} \operatorname{d}\!y \nonumber\\
    &\leq&\frac{-Cm_k(C-\varepsilon_k)}{\operatorname{dist}(\bar{x}^k, T_{\lambda_k})^{2s}}\nonumber\\
    &\leq&\frac{-Cm_k(C-\varepsilon_k)}{2^{2s}}
\end{eqnarray}
On the other hand, combining the differential equation in \eqref{MR2}, \eqref{MR7} with the boundedness of $C_{\lambda_k}(\bar{x}^k,\bar{t}_k)$ and Lemma \ref{lem1}\,, we derive
\begin{eqnarray}\label{MR9}
  (\partial_t-\Delta)^s W_k(\bar{x}^k,\bar{t}_k)&=& (\partial_t-\Delta)^s w_{\lambda_k}(\bar{x}^k,\bar{t}_k)-\varepsilon_km_k(\partial_t-\Delta)^s \eta(\bar{t}_k) \nonumber\\
  &=& C_{\lambda_k}(\bar{x}^k,\bar{t}_k)w_{\lambda_k}(\bar{x}^k,\bar{t}_k)-\varepsilon_km_k\partial_t^s\eta(\bar{t}_k)\nonumber\\
  &\geq&-C_{\lambda_k}(\bar{x}^k,\bar{t}_k)m_k-C\varepsilon_km_k.
\end{eqnarray}
Thus, a combination of \eqref{MR8} and \eqref{MR9} yields that
\begin{equation*}
  C_{\lambda_k}(\bar{x}^k,\bar{t}_k)\geq -C\varepsilon_k+\frac{C(C-\varepsilon_k)}{2^{2s}}\geq C(C-\varepsilon_k).
\end{equation*}
By virtue of $\varepsilon_k\rightarrow 0$ as $k\rightarrow\infty$, we deduce that
\begin{equation*}
C_{\lambda_k}(\bar{x}^k,\bar{t}_k)\geq C_0>0
\end{equation*}
for sufficiently large $k$.
From this, owing to
$$C_{\lambda_k}(\bar{x}^k,\bar{t}_k)=f'(\xi)\,\, \mbox{for}\,\, \xi\in(u_{\lambda_k}(\bar{x}^k,\bar{t}_k),u(\bar{x}^k,\bar{t}_k))$$
and the assumption $f'(0)\leq 0$, there must exist a subsequence of $\{(\bar{x}^k,\bar{t}_k)\}$ (still denoted by $\{(\bar{x}^k,\bar{t}_k)\}$) such that
\begin{equation}\label{MR10}
  u(\bar{x}^k,\bar{t}_k)\geq C_1>0.
\end{equation}

To proceed, we denote
 \begin{equation*}
\bar{w}_{k}(x,t)=w_{\lambda_k}(x,t+\bar{t}_k)\,\,\mbox{and}\,\,\bar{C}_{k}(x,t)= C_{\lambda_k}(x,t+t_k).
\end{equation*}
It follows from
Arzel\`{a}-Ascoli theorem that there exists some functions $\bar {w}(x,t)$ and $\bar C(x,t)$ such that
\begin{equation*}
\lim\limits_{k\rightarrow\infty}\bar{w}_{k}(x,t)=\bar{w}(x,t)\,\, \mbox{and}\,\,\lim\limits_{k\rightarrow\infty}\bar{C}_{k}(x,t)=\bar {C}(x,t).
\end{equation*}
Furthermore,
on account of the equation
\begin{eqnarray*}
(\partial_t-\Delta)^s \bar{w}_{k}(x,t) =\bar{C}_{k}(x,t)\bar{w}_{k}(x,t),\,\, \mbox{in}\,\, \Omega_{\lambda_k}\times \mathbb{R},
\end{eqnarray*}
and
utilizing the regularity theory for master equation established in \cite{ST}, we obtain that the limit function $\bar{w}$ satisfies
\begin{eqnarray}\label{limequ}
(\partial_t-\Delta)^s \bar {w}(x,t) = \bar C(x,t)\bar {w}(x,t),\,\, \mbox{in}\,\, \Omega_{\lambda_0}\times \mathbb{R},
\end{eqnarray}
by $\lambda_k\rightarrow\lambda_0$ as $k\rightarrow\infty$.
Due to $\Omega_{\lambda_k}$ is a bounded domain, we may
assume that $\bar{x}^k\rightarrow x^0$, then applying \eqref{MR7} to derive
\begin{equation*}
 \bar w(x^0,0)\leftarrow \bar{w}_{k}(\bar{x}^k,0)=w_{\lambda_k}(\bar{x}^k,\bar{t}_k)\rightarrow0,
\end{equation*}
as $k\rightarrow\infty$, that is to say,
\begin{equation*}
  \bar w(x^0,0)=0.
\end{equation*}
Combining the limit equation \eqref{limequ} with a direct calculation, we obtain
\begin{eqnarray*}
  0 &=& (\partial_t-\Delta)^s \bar {w}(x^0,0) \\
   &=&-C_{n,s}\int_{-\infty}^{0}\int_{\mathbb{R}^n}
  \frac{\bar {w}(y,\tau)}{(-\tau)^{\frac{n}{2}+1+s}}e^{\frac{|x^0-y|^2}{4\tau}}\operatorname{d}\!y\operatorname{d}\!\tau \\
  &=&C_{n,s}\int_{-\infty}^{0}\int_{\Sigma_{\lambda_0}} \bar {w}(y,\tau)\left[
  \frac{e^{\frac{|x^0-y^{\lambda_0}|^2}{4\tau}}}{(-\tau)^{\frac{n}{2}+1+s}}- \frac{e^{\frac{|x^0-y|^2}{4\tau}}}{(-\tau)^{\frac{n}{2}+1+s}}\right]\operatorname{d}\!y\operatorname{d}\!\tau.
\end{eqnarray*}
Thereby the nonnegativity of $\bar {w}(x,t)$ in $\Sigma_{\lambda_0}\times\mathbb{R}$, the antisymmetry of $\bar {w}(x,t)$ with respect to $x$, and the radial decrease of the kernel ensure that the following identity
\begin{equation}\label{w=0}
\bar w(x,t)\equiv0\,\,\mbox{in}\,\, \mathbb{R}^n\times(-\infty,0]
\end{equation}
holds.

Furthermore, taking the same translation for $u$ as follows
 \begin{equation*}
u_{k}(x,t)=u(x,t+t_k).
\end{equation*}
Similarly to the above argument regarding $\bar{w}_k$, we also have
\begin{equation*}
\lim_{k\rightarrow\infty}u_{k}(x,t)=\bar{u}(x,t),
\end{equation*}
and $\bar{u}(x,t)$ satisfies the limit equation
\begin{eqnarray}\label{limequ2}
(\partial_t-\Delta)^s \bar {u}(x,t) =f(\bar {u}(x,t)),\, {(x, t)\in B_1(0)\times \mathbb{R}}.
\end{eqnarray}
Then it follows from \eqref{MR10} that
\begin{equation}\label{MR11-1}
  \bar{u}(x^0,0)=\lim_{k\rightarrow\infty} u_k(\bar{x}^k,0)=\lim_{k\rightarrow\infty} u(\bar{x}^k,\bar{t}_k)>0.
\end{equation}
Now we claim that
\begin{equation}\label{MR11}
  \bar{u}(x,0)>0\,\, \mbox{for any}\,\, x\in B_1(0).
\end{equation}
If not, then there exists a point $\bar{x}\in B_1(0)$ such that
$$\bar{u}(\bar{x},0)=0=\min_{\mathbb{R}^n\times\mathbb{R}}\bar{u}(x,0)$$
by the exterior condition and the interior positivity of $u$.
It follows that
 \begin{eqnarray*}
(\partial_t-\Delta)^s \bar {u}(\bar{x}, 0)
 =-C_{n,s}\int_{-\infty}^{0}\int_{\mathbb{R}^n}
  \frac{\bar {u}(y,\tau)}{(-\tau)^{\frac{n}{2}+1+s}}e^{\frac{|\bar{x}-y|^2}{4\tau}}\operatorname{d}\!y\operatorname{d}\!\tau \leq0.
\end{eqnarray*}
In contrary, we apply the assumption $f(0)\geqslant0$ and the limit equation \eqref{limequ2} to lead to
\begin{eqnarray*}
(\partial_t-\Delta)^s \bar {u}(\bar{x}, 0)=f(\bar {u}(\bar{x}, 0))=f(0)\geq 0.
\end{eqnarray*}
Thus, we conclude that
\begin{eqnarray*}
(\partial_t-\Delta)^s \bar {u}(\bar{x}, 0)=0.
\end{eqnarray*}
Taking into account that $\bar{u} \geq 0$, we arrive at
\begin{eqnarray*}
\bar {u}(x, t)\equiv0 \,\, \mbox{in}\,\, \mathbb{R}^n\times(-\infty,0],
\end{eqnarray*}
which contradicts \eqref{MR11-1}. From this, we verify that the assertion \eqref{MR11} is valid.

Finally, a combination of $\bar{u}(x,0)\equiv 0$ in $B_1^c(0)$, \eqref{MR11} and $\lambda_0<0$ yields that there must exist $x\in B_1^c(0)$ such that $x^{\lambda_0}\in B_1(0)$ and
$$\bar{w}(x,0)=\bar{u}(x^{\lambda_0},0)-\bar{u}(x,0)=\bar{u}(x^{\lambda_0},0)>0,$$
which contradict \eqref{w=0}. Therefore, we prove that the limiting position must be $T_0$, i.e., $\lambda_0=0$.
By arbitrarily choosing the direction of $x_1$ and combining with the definition of $\lambda_0$, we conclude that $u(x,t)$ must be radially symmetric and monotone decreasing about the origin in $x\in B_1(0)$ for any $t\in\mathbb{R}$.

We have yet to prove that such decrease is strict, in fact, it is sufficient to argue that
\begin{equation}\label{MR12}
  w_\lambda(x,t)> 0 \,\, \mbox{in}\,\,  \Omega_\lambda\times\mathbb{R}
\end{equation}
for any $-1<\lambda<0$. If not, then there exist some $\lambda_0\in(-1,0)$ and a point $(x^0,t_0)\in \Omega_{\lambda_0}\times\mathbb{R}$ such that
\begin{equation*}
  w_{\lambda_0}(x^0,t_0)=0.
\end{equation*}
Combining the differential equation in \eqref{MR2} with the definition of nonlocal operator $(\partial_t-\Delta)^s$, we deduce that
\begin{eqnarray*}
  0 =(\partial_t-\Delta)^s w_{\lambda_0}(x^0,t_0)&=& C_{n,s}\int_{-\infty}^{t_0}\int_{\mathbb{R}^n}
  \frac{-w_{\lambda_0}(y,\tau)}{(t_0-\tau)^{\frac{n}{2}+1+s}}e^{-\frac{|x^0-y|^2}{4(t_0-\tau)}}\operatorname{d}\!y\operatorname{d}\!\tau\\
   &=&  C_{n,s}\int_{-\infty}^{t_0}\int_{\Sigma_{\lambda_0}}
  \frac{w_{\lambda_0}(y,\tau)}{(t_0-\tau)^{\frac{n}{2}+1+s}}\left[e^{-\frac{|x^0-y^{\lambda_0}|^2}{4(t_0-\tau)}}-e^{-\frac{|x^0-y|^2}{4(t_0-\tau)}}\right]\operatorname{d}\!y\operatorname{d}\!\tau.
\end{eqnarray*}
Since $w_{\lambda_0}(x,t)\geq 0$ in $\Sigma_{\lambda_0}\times\mathbb{R}$ and
$$e^{-\frac{|x^0-y^{\lambda_0}|^2}{4(t_0-\tau)}}-e^{-\frac{|x^0-y|^2}{4(t_0-\tau)}}<0,$$
then we must have
$w_{\lambda_0}(x,t)\equiv 0$ in $\Sigma_{\lambda_0}\times(-\infty,t_0]$. However, it contradicts the fact that $w_{\lambda_0}(x,t)\not\equiv 0$ in $\Sigma_{\lambda_0}$ for any fixed $t\in (-\infty,t_0]$, due to the exterior condition and the interior positivity of $u(x,t)$. Hence, we verify that the assertion \eqref{MR12} is valid, and thus the proof of Theorem \ref{Ballsym} is completed.
\end{proof}
\subsection{Liouville type theorem in the whole space}\label{4.2}

At the end of this paper, we utilize the maximum principle in unbounded domains (Theorem \ref{MPUB}) to show that solutions are symmetric with respect to all hyperplanes in $\mathbb{R}^n$ for any $t\in\mathbb{R}$, and
combine with Fourier transform to complete the proof of Liouville theorem (Theorem \ref{Liouville}) for homogeneous master equation
\begin{equation*}
    (\partial_t-\Delta)^s u(x,t)=0 ,\,\,  \mbox{in}\,\,  \mathbb{R}^n\times\mathbb{R}.
\end{equation*}

\begin{proof}
[\bf Proof of Theorem \ref{Liouville}\,.] \,
For any fixed $t\in \mathbb{R}$, we first claim that $u(x,t)$ is symmetric with respect to any hyperplane in $\mathbb{R}^n$.
Let $x_1$ be any given direction in $\mathbb{R}^n$, where we keep the notations $T_\lambda$, $\Sigma_\lambda$, $x^\lambda$, $u_{\lambda}$ and $w_{\lambda}$ as defined above.
For any $\lambda\in \mathbb{R}$, according to equation \eqref{model}, we immediately calculate that $w_\lambda$ satisfies
\begin{equation*}
\left\{
\begin{array}{ll}
    (\partial_t-\Delta)^sw_\lambda(x,t)=0 ,~   & \mbox{in}\,\, \Sigma_\lambda\times\mathbb{R}, \\
    w_\lambda(x,t)=- w_\lambda(x^\lambda,t), &\mbox{in}\,\, \Sigma_\lambda\times\mathbb{R}.
\end{array}
\right.
\end{equation*}
Meanwhile, the boundedness of $u$ implies that $w_\lambda$ is also bounded. Then in terms of Theorem \ref{MPUB}\,, we derive
\begin{equation*}
  w_\lambda(x,t)\leq0 \,\,\mbox{in} \,\,  \Sigma_\lambda\times\mathbb{R}.
\end{equation*}
Replacing $w_\lambda$ with $-w_\lambda$ and following a similar argument as above, we can obtain
\begin{equation*}
  w_\lambda(x,t)\geq0 \,\,\mbox{in} \,\,  \Sigma_\lambda\times\mathbb{R}.
\end{equation*}
Thus, it follows that
\begin{equation*}
  w_\lambda(x,t)=0 \,\,\mbox{in} \,\,  \Sigma_\lambda\times\mathbb{R}.
\end{equation*}
Thereby the arbitrariness of $\lambda$ implies that $u(x,t)$ is symmetric with respect to any hyperplane perpendicular to $x_1$-axis. Furthermore, since the choice of $x_1$ direction is also arbitrary, we verify that $u(x,t)$ is symmetric with respect to any hyperplane in $\mathbb{R}^n$ for any fixed $t\in\mathbb{R}$.
Therefore, we deduce that $u$ must depend only on $t$, i.e.,
\begin{equation*}
  u(x,t)=u(t)\,\,\mbox{in} \,\,  \mathbb{R}^n\times\mathbb{R}.
\end{equation*}

From this, the proof Theorem \ref{Liouville} boils down to showing that the bounded solution $u(t)$ of
\begin{equation}\label{liou3}
D_{\rm left}^s u(t) =0 \,\,\mbox{in}\,\, \mathbb{R}
\end{equation}
must be a constant.

Note that the boundedness of $u$ implies that $u$ belongs to a one-side distributional space ${\mathcal L}^-_{s}(\mathbb{R})$. More precisely,
$$ u(t)\in {\mathcal L}^-_{s}(\mathbb{R})=\{u \in L^1_{\rm loc} (\mathbb{R}) \mid \int_{-\infty}^t \frac{|u(\tau)|}{1+|\tau|^{1+s}}\operatorname{d}\!\tau<+\infty\,\, \mbox{for any} \,\, t\in\mathbb{R}\},$$
in which one can define $D_{\rm left}^s u$ as a distribution
$$\int_{-\infty}^{+\infty}\left(D_{\rm left}^s u(t)\right)\psi(t)\operatorname{d}\!t =\int_{-\infty}^{+\infty} u(t)D_{\rm right}^s\psi(t)\operatorname{d}\!t$$
for any $\psi\in \mathcal{S}$ (cf. \cite{SV}). Here $D_{\rm right}^s$ is the Marchaud right fractional derivative, defined as
\begin{equation*}
D_{\rm right}^s \psi(t)=\frac{1}{|\Gamma(-s)|}\int_{t}^{+\infty} \frac{\psi(t)-\psi(\tau)}{(\tau-t)^{1+s}}\operatorname{d}\!\tau,
\end{equation*}
which only takes into account the values of $\psi$ that occur after time $t$ in the future.
Moreover, in such a setting $u$ is a tempered distribution, then we can define its Fourier transform and the inverse Fourier transform in the sense of distributions.
Applying the fact presented in \cite{SV} that
\begin{equation*}
  \mathcal{F}(D_{\rm right}^s\psi)(\rho)=(-i\rho)^s \mathcal{F}(\psi)(\rho)
\end{equation*}
for any $\psi\in \mathcal{S}$, then it follows from \eqref{liou3} that
\begin{eqnarray}\label{liou1}
   0=\int_{-\infty}^{+\infty}\left(D_{\rm left}^s u(t)\right)\psi(t)\operatorname{d}\!t &=&\int_{-\infty}^{+\infty} u(t)D_{\rm right}^s\psi(t)\operatorname{d}\!t\nonumber \\
  &=& \int_{-\infty}^{+\infty} u(t)\mathcal{F}^{-1}\left((-i\rho)^s \mathcal{F}(\psi)(\rho)\right)(t)\operatorname{d}\!t
\end{eqnarray}
for any $\psi\in \mathcal{S}$.

In the sequel, we show that
\begin{equation}\label{liou2}
  \langle\mathcal{F}u,\phi\rangle=0\,\,\mbox{for any}\,\, \phi\in C_0^\infty(\mathbb{R}\setminus \{0\}).
\end{equation}
Let $\phi\in C_0^\infty(\mathbb{R}\setminus \{0\})$, then the function $\frac{\phi(\rho)}{(-i\rho)^s}$ also belongs to $C_0^\infty(\mathbb{R}\setminus \{0\})\subset\mathcal{S}$. There must exist a function $\psi\in\mathcal{S}$ such that
\begin{equation*}
  \mathcal{F}(\psi)(\rho)=\frac{\phi(\rho)}{(-i\rho)^s}.
\end{equation*}
It follows from \eqref{liou1} that
\begin{eqnarray*}
 \langle\mathcal{F}u,\overline{\phi}\rangle&=& \langle\mathcal{F}u,\overline{(-i\rho)^s\mathcal{F}(\psi)(\rho)}\rangle \\
   &=&  \langle u,\overline{\mathcal{F}^{-1}\left((-i\rho)^s\mathcal{F}(\psi)(\rho)\right)}\rangle\\
  &=& \int_{-\infty}^{+\infty} u(t)\overline{\mathcal{F}^{-1}\left((-i\rho)^s \mathcal{F}(\psi)(\rho)\right)(t)}\operatorname{d}\!t=0.
\end{eqnarray*}
Hence, the assertion \eqref{liou2} is valid, which implies that
$\mathcal{F}(u)$ is supported at the origin. From this, we conclude that $u(t)$ is a polynomial of $t$. While the boundedness of $u$ indicates that $$u(t)\equiv C.$$
In conclusion, we complete the proof of Theorem \ref{Liouville}\,.
\end{proof}

\section*{Acknowledgments}
The authors are grateful for the fruitful discussions with professor Wenxiong Chen (Yeshiva University) during the preparation of this paper. The work of the first author is partially supported by the National Natural Science Foundation of China (NSFC Grant No. 12101452), and the work of the third author is partially supported by the National Natural Science Foundation of China (NSFC Grant No. 12071229).

\section*{Conflict of interest}
The authors declare that they have no conflict of interest.

\end{document}